# About Three Dimensional Double-Sided Dirichlet and Neumann Boundary Value Problems for the Laplacian

Olexandr Polishchuk
Laboratory of Modeling and Optimization of Complex Systems
Pidstryhach Institute for Applied Problems of Mechanics and Mathematics, National Academy of Sciences of Ukraine
Lviv, Ukraine
od_polishchuk@ukr.net

*Abstract*— The orthogonality of Hilbert spaces whose elements can be represented as simple and double layer potentials is determined. Conditions of well-posed solvability of integral equations for the sum of simple and double layer potentials equivalent to double-sided Dirichlet, Neumann, and Dirichlet-Neumann boundary value problems for the Laplacian are established in the Hilbert space, elements of which as well as their normal derivatives have the jump through boundary surface. The properties of boundary operators that relate the double-sided boundary conditions of different types for the three dimensional Laplace equation are investigated.

## I. Introduction

Many physical processes (e.g. diffusion, heat flux, electrostatic field, perfect fluid flow, elastic motion of solid bodies, groundwater flow, etc.) are modeled using boundary value problems for Laplace equation [1, 2]. The powerful tools for solving such problems are potential theory methods, especially in the case of complex shape surface or tired boundary surface [3, 4]. In number of cases, application of potential theory methods requires determination the conditions of well-posed solvability for corresponding integral equations. Review of such conditions for main three dimensional boundary value problems for the Laplacian and equivalent to them integral equations for the simple and double layer potentials contains in [5]. These results allow us to use projection methods [6, 7] for numerical solution of such integral equations, avoiding the use of resource-consuming regularization procedures [8, 9]. The need to determine the conditions of well-posed solvability also arises when the sum of simple and double layer potentials is used to solve the double-sided Dirichlet and Neumann problems or double-sided Dirichlet-Neumann problem in the space of functions that, same as their normal derivatives, have a jump on crossing boundary surface. These conditions allow us to construct and determine the properties of so-called boundary operators, which put in accordance known values of the boundary condition of one type the values of the boundary condition of another type [10]. Obtained results enable to build effective methods for numerical solution of investigated problems.

## II. Functional spaces $H^1_{\Gamma,\Delta=0}$, $K^1_{\Gamma,\Delta=0}$, and $HK^1_{\Gamma,\Delta=0}$

Let $G$ be the bounded open set in $R^3$, the boundary of which is Lipshitz surface $\Gamma$. Let us denote $G' = R^3 \setminus \overline{G}$ and introduce into $G$, $G'$ Sobolev spaces $W_2^1(G)$ and

$$W_{2,0}^1(G') = \{u^e \in D'(G') : u^e/r, Du^e \in L_2(G')\}.$$

Let us determine the Hilbert space $W_2^{1/2}(\Gamma)$ on the surface $\Gamma$ and introduce the space

$$H^1 = W_2^1(G) \times W_{2,0}^1(G').$$

Determine the linear continuous trace operators [11]

$$\gamma_0^i u^i = u^i\big|_{\Gamma_i},\ \gamma_0^e u^e = u^e\big|_{\Gamma_e},\ \gamma_0^i : W_2^1(G) \to W_2^{1/2}(\Gamma),\ \gamma_0^e : W_{2,0}^1(G') \to W_2^{1/2}(\Gamma),$$

$$\gamma_1^i u^i = \partial u^i/\partial \mathbf{n}\big|_{\Gamma_i},\ \gamma_1^e u^e = \partial u^e/\partial \mathbf{n}\big|_{\Gamma_e},\ \gamma_1^i : W_2^1(G) \to W_2^{-1/2}(\Gamma),\ \gamma_1^e : W_{2,0}^1(G') \to W_2^{-1/2}(\Gamma)$$

where $\Gamma_i$ and $\Gamma_e$ are the internal and external sides of surface $\Gamma$, accordingly, $\mathbf{n}$ is a normal to the surface $\Gamma$ external to the domain $G$, and $W_2^{-1/2}(\Gamma)$ is the space duel to $W_2^{1/2}(\Gamma)$. Denote

$$\gamma_0^1 u = \gamma_0^i u^i - \gamma_0^e u^e,\ \gamma_1^1 u = \gamma_1^i u^i - \gamma_1^e u^e,$$

where $u = (u^i, u^e)$, $u^i \in W_2^1(G)$, $u^e \in W_{2,0}^1(G')$, and

$$W_2^1(G; \Delta = 0) = \{u^i \in W_2^1(G) : \Delta u^i = 0\},$$

$$W_{2,0}^1(G'; \Delta = 0) = \{u^e \in W_{2,0}^1(G') : \Delta u^e = 0\}.$$

For arbitrary $u^i \in W_2^1(G; \Delta = 0)$ and $v^i \in W_2^1(G)$ we have the Green's formula [12]

$$\int_G \nabla u^i \nabla v^i = \left\langle \gamma_1^i u^i, \gamma_0^i v^i \right\rangle_{W_2^{-1/2}(\Gamma) \times W_2^{1/2}(\Gamma)} \tag{1}$$

where $\langle .,. \rangle_{V \times V^*}$ is the duality relation on $V \times V^*$. For arbitrary $u^e \in W_{2,0}^1(G'; \Delta u = 0)$ and $v^e \in W_{2,0}^1(G')$ we have the Green's formula [12]

$$\int_{G'} \nabla u^e \nabla v^e = -\left\langle \gamma_1^e u^e, \gamma_0^e v^e \right\rangle_{W_2^{-1/2}(\Gamma) \times W_2^{1/2}(\Gamma)}. \tag{2}$$

Introduce the space

$$H_{\Delta=0}^1 = W_2^1(G; \Delta = 0) \times W_{2,0}^1(G'; \Delta u = 0).$$

Then from (1) and (2) we obtain the Green's formula

$$\int_{G \cup G'} \nabla u \nabla v = \left\langle \gamma_1^i u, \gamma_0^i v \right\rangle_{W_2^{-1/2}(\Gamma) \times W_2^{1/2}(\Gamma)} - \left\langle \gamma_1^e u, \gamma_0^e v \right\rangle_{W_2^{-1/2}(\Gamma) \times W_2^{1/2}(\Gamma)} \tag{3}$$

for arbitrary $u \in H_{\Delta=0}^1$ and $v \in H^1$.

Introduce the Hilbert space

$$H_\Gamma^1 = \{u \in H^1 : \gamma_0^1 u = 0\}, \quad (u,v)_{H_\Gamma^1} = \int_{G \cup G'} \nabla u \nabla v, \quad \|u\|_{H_\Gamma^1} = (u,u)_{H_\Gamma^1}^{1/2},$$

where $(.,.)_V$ is the scalar product on $V$. Introduce the space

$$H_{\Gamma, \Delta=0}^1 = \{u \in H_\Gamma^1 : \Delta u = 0\}.$$

For arbitrary $u \in H_{\Gamma, \Delta=0}^1$ and $v \in H_\Gamma^1$ from (1) and (2), we have the Green's formula

$$(u,v)_{H_\Gamma^1} = \left\langle \gamma_1^1 u, \gamma_0 v \right\rangle_{W_2^{-1/2}(\Gamma) \times W_2^{1/2}(\Gamma)} \tag{4}$$

where $\gamma_0 v = \gamma_0^i v^i = \gamma_0^e v^e$.

The next result is in order [12].

Theorem 1. Operator $\gamma_1^1$ is an isomorphism of $H_{\Gamma, \Delta=0}^1$ onto $W_2^{-1/2}(\Gamma)$, and arbitrary function $u \in H_{\Gamma, \Delta=0}^1$ can be represented as

$$u(x) = (U\gamma_1^1 u)(x) \equiv \frac{1}{4\pi} \int_\Gamma \frac{\gamma_1^1 u(y)}{|x-y|} d\Gamma_y, \quad x \in G, G', \quad y \in \Gamma.$$

Introduce the Hilbert space [13]

$$K_\Gamma^1 = \{u \in H^1 \setminus R : \gamma_1^1 u = 0\}$$

where $R$ is the space of functions that are constant on $G$ and vanishing on $G'$, and

$$(u,v)_{K_\Gamma^1} = \int_{G \cup G'} \nabla u \nabla v, \quad \|u\|_{K_\Gamma^1} = (u,u)_{K_\Gamma^1}^{1/2}.$$

Introduce the spaces

$$K_{\Gamma, \Delta=0}^1 = \{u \in K_\Gamma^1 : \Delta u = 0\},$$

$$\widetilde{W}_2^{1/2}(\Gamma) = W_2^{1/2}(\Gamma) \setminus P,$$

$$\widetilde{W}_2^{-1/2}(\Gamma) = \{g \in W_2^{-1/2}(\Gamma) : \langle g, 1 \rangle_{W_2^{1/2}(\Gamma) \times W_2^{-1/2}(\Gamma)} = 0\}$$

where $P$ is the set of constant functions on $\Gamma$. For arbitrary $u \in K_{\Gamma,\Delta=0}^1$ and $v \in K_\Gamma^1$ from (1) and (2), we have the Green's formula

$$(u,v)_{K_\Gamma^1} = \langle \gamma_0^1 u, \gamma_1 v \rangle_{W_2^{1/2}(\Gamma) \times W_2^{-1/2}(\Gamma)} \tag{5}$$

where $\gamma_1 v = \gamma_1^i v^i = \gamma_1^e v^e$.

The next result is in order [13].

Theorem 2. Operator $\gamma_0^1$ is an isomorphism of $K_{\Gamma,\Delta=0}^1$ onto $\widetilde{W}_2^{1/2}(\Gamma)$, and arbitrary function $v \in K_{\Gamma,\Delta=0}^1$ can be represented as

$$v(x) = (V \gamma_0^1 v)(x) \equiv -\frac{1}{4\pi} \int_\Gamma \gamma_0^1 v(y) \frac{\partial}{\partial \mathbf{n}_y} \frac{1}{|x-y|} d\Gamma_y, \ x \in G, G', \ y \in \Gamma.$$

Assume that

$$u \in H_{\Gamma,\Delta=0}^1 \cap K_{\Gamma,\Delta=0}^1.$$

Since $u \in H_{\Gamma,\Delta=0}^1$, i.e. $\gamma_0^1 u = 0$, from (3) we have

$$\int_{G \cup G'} (\nabla u)^2 = \langle \gamma_1^1 u, \gamma_0 u \rangle_{W_2^{-1/2}(\Gamma) \times W_2^{1/2}(\Gamma)}.$$

But $u \in K_{\Gamma,\Delta=0}^1$ also, i.e. $\gamma_1^1 u = 0$. Thus,

$$\int_{G \cup G'} (\nabla u)^2 = 0.$$

Like that the spaces $H_{\Gamma,\Delta=0}^1$ and $K_{\Gamma,\Delta=0}^1$ are ortogonal relatively introducing scalar product. Introduce the spaces

$$HK_\Gamma^1 = H_\Gamma^1 \oplus K_\Gamma^1, \ (u,v)_{HK_\Gamma^1} = \int_{G \cup G'} \nabla u \nabla v, \ \|u\|_{HK_\Gamma^1} = (u,u)_{HK_\Gamma^1}^{1/2}$$

and

$$HK_{\Gamma,\Delta=0}^1 = H_{\Gamma,\Delta=0}^1 \oplus K_{\Gamma,\Delta=0}^1.$$

Then from theorems 1 and 2 we have that the following result is in order [14].

Theorem 3. Operator $(\gamma_1^1, \gamma_0^1)$ is an isomorphism of $HK_{\Gamma,\Delta=0}^1$ onto $W_2^{-1/2}(\Gamma) \times \widetilde{W}_2^{1/2}(\Gamma)$ and arbitrary function $w \in HK_{\Gamma,\Delta=0}^1$ can be singlevaluedly represent as

$$w(x) = (W(\gamma_1^1 u, \gamma_0^1 v))(x) = (U, V) \ (\gamma_1^1 u, \gamma_0^1 v)^T (x) = (U \ \gamma_1^1 u)(x) + (V \ \gamma_0^1 v)(x) = u(x) + v(x), \tag{6}$$

where

$$u \in H_{\Gamma,\Delta=0}^1 \tag{7}$$

and

$$v \in K_{\Gamma,\Delta=0}^1, \ x \in G, G'. \tag{8}$$

III. DIRICHLET AND NEUMANN BOUNDARY VALUE PROBLEMS FOR THE LAPLACIAN IN $R^3$

Consider the next boudary value problems for the Laplacian .

1. Double-sided Dirichlet problem: to find the function

$$w = (w^i, u^e) \in HK_{w,\Delta=0}^1, \tag{9}$$

which satisfies conditions

$$\gamma_0^i w^i = f^i, \ \gamma_0^e w^e = f^e, \ f^i, f^e \in W_2^{1/2}(\Gamma). \tag{10}$$

Denote

$$\gamma_0^1 w = q = f^i - f^e. \qquad (11)$$

Account that an arbitrary element $w \in HK^1_{\Gamma, \Delta=0}$ can be singlevaluedly represents as a sum (6) we obtain from (11)

$$\gamma_0^1 v = q. \qquad (12)$$

From theorem 2 we have that following result is in order

Theorem 4. If

$$q \in \widetilde{W}_2^{1/2}(\Gamma), \qquad (13)$$

the problem (8), (12) has one and only one solution.

Let us the condition (13) is carrying out and the values $\gamma_0^i v^i = \varphi^i$ or $\gamma_0^e v^e = \varphi^e$ are defined. Then

$$\gamma_0 u = f^e - \varphi^e = f^i - \varphi^i \in W_2^{1/2}(\Gamma) \qquad (14)$$

and from theorem 1 we have that following result is in order.

Theorem 5. The problem (7), (14) has one and only one solution.

Like that we have

Theorem 6. If

$$q = f^i - f^e \in \widetilde{W}_2^{1/2}(\Gamma),$$

the problem (9)-(10) has one and only one solution.

Consider in $HK^1_{\Gamma, \Delta=0}$ an integral equation for the sum of simple and double layer potentials equivalents to the double-sided Dirichlet problem for the Laplacian. Reduce the problem (9), (10) to an integral equation using the first of equality (10), equality (6), result of the theorem 3 and formula

$$(\gamma_0^i v)(x) = -q(x)/2 - \frac{1}{4\pi} \int_\Gamma q(y) \frac{\partial}{\partial \vec{n}_y} \frac{1}{|x-y|} d\Gamma_y, \quad x \in \Gamma$$

Then we have

$$\gamma_0^i w^i(x) = \gamma_0 u(x) + \gamma_0^i v^i(x) = \frac{1}{4\pi} \int_\Gamma \frac{\sigma(y)}{|x-y|} d\Gamma_y - q(x)/2 - \frac{1}{4\pi} \int_\Gamma q(y) \frac{\partial}{\partial \vec{n}} \frac{1}{|x-y|} d\Gamma_y = f^i(x), \quad x \in \Gamma$$

So far as the value $q(x)$ is known we come to integral equation of the first kind

$$\frac{1}{4\pi} \int_\Gamma \frac{\sigma(y)}{|x-y|} d\Gamma_y = f^i(x) + q(x)/2 + \frac{1}{4\pi} \int_\Gamma q(y) \frac{\partial}{\partial \vec{n}} \frac{1}{|x-y|} d\Gamma_y, \quad x \in \Gamma \qquad (15)$$

Then analogously to [13, 15] we proves

Theorem 7. The equation (15) has one and only one solution.

2. Double-sided Neumann problem: to find the function $w \in HK^1_{\Gamma, \Delta=0}$ which satisfies conditions

$$\gamma_1^i w^i = g^i, \ \gamma_1^e w^e = g^e, \ g^i, g^e \in W_2^{-1/2}(\Gamma) \qquad (16)$$

Denote

$$\gamma_1^1 w = \sigma = g^i - g^e \qquad (17)$$

Account that an arbitrary element $w \in HK^1_{\Gamma, \Delta=0}$ can be represent as (6)-(8) we obtain from (17)

$$\gamma_1^1 u = \sigma \in W_2^{-1/2}(\Gamma) \qquad (18)$$

From theorem 1 we have that following result is in order.

Theorem 8. The problem (7), (18) has one and only one solution.

Let the values $\gamma_1^i u^i = \phi^i$ or $\gamma_1^e u^e = \phi^e$ are defined. Then

$$\gamma_1 v = g^i - \phi^i = g^e - \phi^e. \qquad (19)$$

From condition
$$H^1_{\Gamma,\Delta=0} \cap K^1_{\Gamma,\Delta=0} = \{0\}$$

we have that the right part of (19) belong to $\widetilde{W}_2^{-1/2}(\Gamma)$. Then we have [13]

Theorem 9. The problem (7), (19) has one and only one solution.

Like that we have

Theorem 10. The problem (9), (16) has one and only one solution.

Reduce the problem (9), (16) to an integral equation using the second of equality (16), equality (6), result of the theorem 9 and formula [13]

$$(\gamma_0^e u)(x) = -\sigma(x)/2 + \frac{1}{4\pi}\int_\Gamma \sigma(y)\frac{\partial}{\partial \vec{n}_y}\frac{1}{|x-y|}d\Gamma_y, \quad x \in \Gamma.$$

Then

$$\gamma_1^e w^e(x) = \gamma_1 v(x) + \gamma_1^e u^e(x) = \frac{1}{4\pi}\frac{\partial}{\partial \vec{n}_x}\int_\Gamma q(y)\frac{\partial}{\partial \vec{n}}\frac{1}{|x-y|}d\Gamma_y - \sigma(x)/2 + \frac{1}{4\pi}\frac{\partial}{\partial \vec{n}_x}\int_\Gamma \frac{\sigma(y)}{|x-y|}d\Gamma_y = g^e(x), \quad x \in \Gamma.$$

So far as the value $\sigma(x)$ is known we come to integral equation of the first kind

$$\frac{1}{4\pi}\frac{\partial}{\partial \vec{n}_x}\int_\Gamma q(y)\frac{\partial}{\partial \vec{n}}\frac{1}{|x-y|}d\Gamma_y = g^e(x) + \sigma(x)/2 - \frac{1}{4\pi}\frac{\partial}{\partial \vec{n}_x}\int_\Gamma \frac{\sigma(y)}{|x-y|}d\Gamma_y, \quad x \in \Gamma \quad (20)$$

Then analogously to [13, 15] we proves

Theorem 11. The equation (20) has one and only one solution.

3. Double-sided Dirichlet-Neumann problems: to find the function $w \in HK^1_{\Gamma,\Delta=0}$ which satisfies conditions

$$\gamma_0^i w^i = f^i, \quad f^i \in W_2^{1/2}(\Gamma), \quad \gamma_1^e w^e = g^e, \quad g^e \in W_2^{-1/2}(\Gamma), \quad (21)$$

or conditions

$$\gamma_1^i w^i = g^i, \quad g^i \in W_2^{-1/2}(\Gamma), \quad \gamma_0^e w^e = f^e, \quad f^e \in W_2^{1/2}(\Gamma). \quad (22)$$

Using formula (3) and conditions (21) we have the equivalent to (9), (21) variational problem: to find the function $w \in HK^1_{\Gamma,\Delta=0}$ which satisfies the ratio

$$a(w,z) \equiv \int_{G \cup G'} \nabla w \nabla z = \left\langle f^i, \gamma_1^i z \right\rangle_{W_2^{1/2}(\Gamma) \times W_2^{-1/2}(\Gamma)} - \left\langle g^e, \gamma_0^e z \right\rangle_{W_2^{-1/2}(\Gamma) \times W_2^{1/2}(\Gamma)}, \quad (23)$$

for arbitrary $z \in HK^1_{\Gamma,\Delta=0}$. Then from Lax-Milgram theorem [16] we obtain [17]

Theorem 12. The problem (9), (21) has one and only one solution.

Similarly, the problem (9), (22) can match the equivalent variational problem: to find the function $w \in HK^1_{\Gamma,\Delta=0}$ which satisfies the ratio

$$a(w,z) = \left\langle g^i, \gamma_0^i z \right\rangle_{W_2^{-1/2}(\Gamma) \times W_2^{1/2}(\Gamma)} - \left\langle f^e, \gamma_1^e z \right\rangle_{W_2^{1/2}(\Gamma) \times W_2^{-1/2}(\Gamma)}, \quad (24)$$

for arbitrary $z \in HK^1_{\Gamma,\Delta=0}$. Then from Lax-Milgram theorem we obtain

Theorem 13. The problem (9), (22) has one and only one solution.

Consider the variational problem (9), (23). Using formula (3) and theorem 3 we have the equivalent to (9), (23) variational problem: to find the pair

$$(\sigma, q) \in W_2^{-1/2}(\Gamma) \times \widetilde{W}_2^{1/2}(\Gamma), \quad \sigma = \gamma_1^1 w, \quad q = \gamma_0^1 w, \quad (25)$$

which satisfies the ratio

$$b((\sigma,q),(\mu,v)) \equiv \int_\Gamma [\frac{1}{4\pi}\int_\Gamma \frac{\sigma(y)}{|x-y|}d\Gamma_y - q(x)/2 - \frac{1}{4\pi}\int_\Gamma q(y)\frac{\partial}{\partial \vec{n}_y}\frac{1}{|x-y|}d\Gamma_y] \times [\mu(x)/2 +$$

$$+\frac{1}{4\pi}\frac{\partial}{\partial \vec{n}_x}\int_\Gamma \frac{\mu(y)}{|x-y|}d\Gamma_y+\frac{1}{4\pi}\frac{\partial}{\partial \vec{n}_x}\int_\Gamma v(y)\frac{\partial}{\partial \vec{n}_y}\frac{1}{|x-y|}d\Gamma_y\,]d\Gamma_x-\int_\Gamma[-\sigma(x)/2+\frac{1}{4\pi}\frac{\partial}{\partial \vec{n}_x}\int_\Gamma \frac{\sigma(y)}{|x-y|}d\Gamma_y+$$

$$+\frac{1}{4\pi}\frac{\partial}{\partial \vec{n}_x}\int_\Gamma q(y)\frac{\partial}{\partial \vec{n}_y}\frac{1}{|x-y|}d\Gamma_y\,]\times[\frac{1}{4\pi}\int_\Gamma \frac{\mu(y)}{|x-y|}d\Gamma_y+v(x)/2-\frac{1}{4\pi}\int_\Gamma v(y)\frac{\partial}{\partial \vec{n}_y}\frac{1}{|x-y|}d\Gamma_y\,]d\Gamma_x= \quad (26)$$

$$=\int_\Gamma f^i(x)[\,\mu(x)/2+\frac{1}{4\pi}\frac{\partial}{\partial \vec{n}_x}\int_\Gamma \frac{\mu(y)}{|x-y|}d\Gamma_y+\frac{1}{4\pi}\frac{\partial}{\partial \vec{n}_x}\int_\Gamma v(y)\frac{\partial}{\partial \vec{n}_y}\frac{1}{|x-y|}d\Gamma_y\,]d\Gamma_x-$$

$$-\int_\Gamma g^e(x)[\,\frac{1}{4\pi}\int_\Gamma \frac{\mu(y)}{|x-y|}d\Gamma_y+v(x)/2-\frac{1}{4\pi}\int_\Gamma v(y)\frac{\partial}{\partial \vec{n}_y}\frac{1}{|x-y|}d\Gamma_y\,]d\Gamma_x,$$

for arbitrary pair $(\mu,v)\in W_2^{-1/2}(\Gamma)\times \widetilde{W}_2^{1/2}(\Gamma)$.

It is obvious that the bilinear form $b((\sigma,q),(\mu,v))$ is continuous and symmetric. Let us $\mu=\sigma$ and $v=q$. Then, using formulas (3), (6) and theorems 1, 2 we obtain

$$b((\sigma,q),(\sigma,q))=\left\langle \gamma_1^i w^i,\gamma_0^i w^i\right\rangle_{W_2^{-1/2}(\Gamma)\times W_2^{1/2}(\Gamma)}-\left\langle \gamma_1^e w^e,\gamma_0^e w^e\right\rangle_{W_2^{-1/2}(\Gamma)\times W_2^{1/2}(\Gamma)}=$$

$$=\|w\|^2_{HK^1_\Gamma}=\|u\|^2_{H^1_\Gamma}+\|v\|^2_{K^1_\Gamma}\geq C_1\|\sigma\|^2_{W_2^{-1/2}(\Gamma)}+C_2\|q\|^2_{\widetilde{W}_2^{1/2}(\Gamma)},\ C_1,C_2\geq 0,$$

for arbitrary pair $(\sigma,q)\in W_2^{-1/2}(\Gamma)\times \widetilde{W}_2^{1/2}(\Gamma)$, i. e. bilinear form $b((\sigma,q),(\mu,v))$ is $W_2^{-1/2}(\Gamma)\times \widetilde{W}_2^{1/2}(\Gamma)$- elliptic.

Obviously, the solution of problem (25)-(26) is equivalent to the solution of the system of integral equations

$$\frac{1}{4\pi}\int_\Gamma \frac{\sigma(y)}{|x-y|}d\Gamma_y-q(x)/2-\frac{1}{4\pi}\int_\Gamma q(y)\frac{\partial}{\partial \vec{n}_y}\frac{1}{|x-y|}d\Gamma_y=f^i(x),\ x\in\Gamma \quad (27)$$

$$-\sigma(x)/2+\frac{1}{4\pi}\frac{\partial}{\partial \vec{n}_x}\int_\Gamma \frac{\sigma(y)}{|x-y|}d\Gamma_y+\frac{1}{4\pi}\frac{\partial}{\partial \vec{n}_x}\int_\Gamma q(y)\frac{\partial}{\partial \vec{n}_y}\frac{1}{|x-y|}d\Gamma_y=g^e(x),\ x\in\Gamma. \quad (28)$$

Then, from the Lax-Milgram theorem, the following result is in order.

Theorem 14. The system (27)-(28) has one and only one solution.

Similarly, solving the problem (9), (26) is equivalent to solving the system of integral equations

$$\sigma(x)/2+\frac{1}{4\pi}\frac{\partial}{\partial \vec{n}_x}\int_\Gamma \frac{\sigma(y)}{|x-y|}d\Gamma_y+\frac{1}{4\pi}\frac{\partial}{\partial \vec{n}_x}\int_\Gamma q(y)\frac{\partial}{\partial \vec{n}_y}\frac{1}{|x-y|}d\Gamma_y=g^i(x),\ x\in\Gamma, \quad (29)$$

$$\frac{1}{4\pi}\int_\Gamma \frac{\sigma(y)}{|x-y|}d\Gamma_y+q(x)/2-\frac{1}{4\pi}\int_\Gamma q(y)\frac{\partial}{\partial \vec{n}_y}\frac{1}{|x-y|}d\Gamma_y=f^e(x),\ x\in\Gamma, \quad (30)$$

and we have

Theorem 15. The system (29)-(30) has one and only one solution.

By means of approximation the spaces $W_2^{-1/2}(\Gamma)$ and $\widetilde{W}_2^{1/2}(\Gamma)$ by complete orthonormal systems of functions, B-splines, or finite elements of other types [6], similarly [7], we can prove the convergence of a number of projection methods (collocation, Galerkin, least squares, least errors) for solution the systems of integral equations (27)-(28) and (29)-(30).

IV. BOUNDARY OPERATORS FOR DOUBLE-SIDED BOUNDARY VALUE PROBLEMS FOR THE LAPLACIAN IN $R^3$

A boundary operator is an operator that allows one to find unknown values of one type of boundary condition by known values of the boundary condition of another type [10]. For the double-sided Dirichlet and Neumann problems, this operator is defined by the formula

$$F(\gamma_1^i w,\gamma_1^e w)=(f^i,f^e) \quad (31)$$

where $(f^i,f^e)$ is the known values of double-sided Dirichlet conditions and $(\gamma_1^i w,\gamma_1^e w)$ is an unknown values of double-sided Neumann conditions. It follows from Theorems 6, 10, 12 and 13 that for all double-sided boundary value problems discussed above, the corresponding boundary operators exist and are isomorphisms on the corresponding pairs of spaces (operator $F$ in formula (31) is an isomorphism from $\widetilde{W}_2^{-1/2}(\Gamma) \times \widetilde{W}_2^{-1/2}(\Gamma)$ onto $W_2^{1/2}(\Gamma) \times W_2^{1/2}(\Gamma)$ etc).

It is impractical to construct boundary operators numerically in the differential formulation, since the use of finite differences or finite elements methods requires the calculation of searched solution in many points of sufficiently large three dimensional domains. Using the integral equations method makes it much easier to find the required values. Indeed, if we know the values of Dirichlet conditions on the boundary surface, then by calculating by formula (11) the values of double layer potential density and solving equation (15), we can calculate the required values of the Neumann conditions on the boundary surface by formula (6). The construction of boundary operators in the integral formulation is correct, since it follows from theorems 7, 11, 14 and 15 that they are isomorphisms on the corresponding pairs of spaces.

Boundary operators are useful in investigating the conditions of well-posed solvability of many problems in mathematical physics [18, 19], for constructing effective methods for their numerical solution [20], and for solving practically important applied problems. Thus, when analyzing the quality of prosthetics of the person 's lower extremities, it is important to determine on the basis of the temperature on the stump's surface the direction of temperature flows (inside or outside) [21]. The answer to this question is given by the construction of boundary operator that puts in accordance to the known Dirichlet boundary value the Neumann boundary value.